\documentstyle[12pt,a4]{article}

\newcommand{\be}{\begin{equation}}
\newcommand{\ee}{\end{equation}}
\newcommand{\bea}{\begin{eqnarray}}
\newcommand{\eea}{\end{eqnarray}}

\newcommand{\ov}{\overline}
\title{Quommutator deformations of osp(2,2)} 
\author{Y. Brihaye \\
Department of Mathematical Physics\\
University of Mons\\
Pl. du Parc, B-7000 MONS, Belgium}
\begin{document}
\begin{titlepage}
\maketitle
\thispagestyle{empty}
\begin{abstract}
We analyse the Witten-Woronowicz's type 
deformations of the Lie superalgebra osp(2,2)
and obtain a deformation parametrized by three
independent parameters. 
For some of these algebras, finite dimensional 
representations are formulated in terms of
a finite difference operator, providing 
operators that are relevant for the classification
of quasi exactly solvable, finite difference equations.
Similar representations are also pointed out for the
Lie superalgebra osp(1,2). 
\end{abstract}
\end{titlepage}
\section{Introduction}
The study of deformations of classical Lie algebras
and superalgebra has attracted much attention in the last few years.
Many types of deformations have been proposed. In some of them,
the (anti)commutators of the generators are set equal to non linear 
(in fact transcendental)
functions of some of the generators, they are the 
Drinfeld-Jimbo deformations \cite{dri,jim}. Another type of deformation
consists in keeping the right hand sides of the structure relations
linear in the generators while the (anti)commutators in the left
hand side are deformed into quommutators.
For the deformations of the second type, also called 
of Witten-Woronowicz's type \cite{wit,wor}, 
all the generators are treated on the same footing
and the quommutator relations provide a set of simple rules to write
the elements of the enveloping algebra in a canonical form; they 
define a normal order. 
This last property is crucial for the kind of applications that we have
in mind, namely the classification of the quasi exactly solvable (QES),
finite difference operators. 
\par The QES (differential or finite difference) 
equations refer to a class of spectral 
equations for which  a finite part of the 
spectrum can be obtained by solving a system
of algebraic equations \cite{tur1,ush}. The QES operators 
defining the QES equations are 
therefore intimately related to the linear 
operators preserving a finite dimensional vector 
space of smooth functions.
Typically this vector space is a direct sum of 
spaces of polynomials of given degree.  
\par There is a close relation between  QES differential operators 
and the  theory of algebras.  Indeed, some  
basic QES operators correspond to the generators,
in a suitable representation, of an abstract
algebra (e.g. in the case of scalar QES operators of one variable, 
it is the Lie algebra sl(2) \cite{tur2}) and the generic 
QES operators appear as the elements of the 
enveloping algebra. 
With the aim to classify the QES operators it is therefore desirable to
possess a set of normal ordering rules for the generators. 
\par If we want to describe the algebraic structures underlying
the QES finite difference equations, the representations
of some deformed algebra seem to emerge in a natural way.
Fo example, the scalar QES  finite difference 
operators in one variable 
\cite{tur2} are related to the Witten deformation of sl(2). 
\par The crucial role played by
the Lie superalgebra osp(2,2) in the study of QES
systems of two equations \cite{turshif,bk} encouraged
us to study the Witten-Woronowicz's
type of deformations of this algebra. This is the topic
of the next section. 
\section{Deformations of osp(2,2)}
The Lie superalgebra osp(2,2) has eight generators,
four of them (the bosonic ones) assemble into a gl(2) 
subalgebra. In order to make easy the comparaison of
our result with ref.\cite{faza} we note these generators
$E_{11}, E_{22},E_{12},E_{21}$. The remaining
(fermionic) generators split into two doublets
under the adjoint action of the gl(2) subalgebra,
we denote them $V_1, V_2$ and $\overline V_1,\overline V_2$. 
\par For later convenience, we define
the quommutator and antiquommutator respectively as
\be
   [A , B]_q = AB - q BA \ \ \ \ , \ \ \ \ \{A , B \}_q = AB + qBA
\label{quomm}
\ee
where $q$ is the deformation parameter.
\par We studied the quommutators deformations of osp(2,2) subject
to the following requirements:
\begin{itemize}
\item each commutator (resp. anticommutator) defining the 
classical algebra is replaced by a quommutator (resp. antiquommutator)
defined in (\ref{quomm}) with its own parameter $q$.
\item the couples of generators which (anti)commute are imposed to
(anti)quommute.
\item the structure of the deformed algebra relates the (anti)quomutators 
of any couple of generators to linear combinations of the generators
but   $[E_{12},E_{21}]_q$ which is allowed to depend
quadratically on the fermionic generators.
\end{itemize}
\par The most general deformation of osp(2,2) obeying these restrictions
and compatible with associativity depends on three parameters, say
$p,r,s$. The different relations read as follows 
\begin{itemize}
\item for the fermionic-fermionic relations
\be
    V_1^2 = V_2^2 = 0 \ \ \ , \ \ \ 
   \overline V_1^2 = \overline V_2^2 = 0
\ee
\be
      \{ V_1,V_2 \}_{{p \over sr}} = 0 \ \ \ , \ \ \ 
      \{ \overline V_1, \overline V_2 \}_{{s \over pr}} = 0 
\ee
\bea
 & \{\overline V_1 , V_1 \} = E_{11} \ \ \ \ , \ \ \ 
 & \{\overline V_2 , V_2 \} = E_{22} \nonumber \\
 & \{\overline V_1 , V_2 \}_{psr} = E_{12} \ \ \ \ , \ \ \ 
 & \{\overline V_2 , V_2 \}_{ps/r} = E_{22}
\label{bff}
\eea
The last equations  can be seen as  defining the 
bosonic operators $E_{ij}$.
\item for the fermionic-bosonic relations
\bea
 \left[ E_{11},V_1 \right]   = 0 & \ \ \   ,
 \ \ \  &  \left[ E_{22},V_1 \right]_{s^2}  = V_1  \nonumber \\
 \left[ E_{21},V_1 \right]_{ps/r}   = 0 & \ \ \ ,
 \ \ \ &  \left[ E_{12},V_1 \right]_{sr/p}  = -{sr\over p} V_2
 \nonumber  \\
  \left[ E_{11},V_2 \right]_{p^2}   = V_2 &  \ \ \ ,
 \ \ \ &  \left[ E_{22},V_2 \right]  = 0  \nonumber \\
  \left[ E_{21},V_2 \right]_{p/sr}  = -{p\over sr} V_1  & \ \ \ , 
   \ \ \ &  \left[ E_{12},V_2 \right]_{psr} = 0 
\eea
\bea
  \left[ E_{11},\ov V_1 \right]  = 0  & \ \ \  , \ \ \
& \left[ E_{22},\ov V_1 \right]_{1/s^2}  = -{1\over s^2} \ov V_1  \nonumber \\
  \left[ E_{21},\ov V_1 \right]_{pr/s}  = \ov V_2  & \ \ \ ,  \ \ \
& \left[ E_{12}, \ov V_1 \right]_{1/psr}  = 0  \nonumber \\
  \left[ E_{11},\ov V_2 \right]_{1/p^2}  = - {1 \over p^2} \ov V_2 & \ \ \  , \ \ \
& \left[ E_{22},\ov V_2 \right]  = 0  \nonumber \\
  \left[ E_{21},\ov V_2 \right]_{r/ps}  = 0 &  \ \ \ , \ \ \
& \left[ E_{12}, \ov V_2 \right]_{s/pr}  = \ov V_1 
\eea
\item and finally for the bosonic-bosonic operators
\be
    [E_{11}, E_{22}] = 0 \ \ , \label{b01} \\
\ee
\bea
  \left[ E_{11}, E_{21} \right]_{1/p^2} = - {1 \over p^2} E_{21} & \ \ \  , \ \ \   
& \left[ E_{22}, E_{21} \right]_{s^2} =  E_{21} \label{b0pm} \\   
  \left[ E_{11}, E_{12} \right]_{p^2} =  E_{12}     &  \  \ \ , \ \ \  
& \left[ E_{22}, E_{12} \right]_{1/s^2} = - {1 \over s^2}  E_{12} 
\label{b1pm}   
\eea
\be
    [E_{22}, E_{12}]_{s^2/p^2} = E_{11} - {s^2 \over p^2} E_{22}
 +(s^2 -1)V_1 \overline V_1 - {s^2 \over p^2}(p^2-1) V_2 \overline V_2
\label{bpm}
\ee
\end{itemize}
The three deformation parameters $p,r,s$ are 
intrinsic and cannot be eliminated by a rescaling of the generators.
This is seen easily from the way these parameters enter in 
the different coefficients defining the quommutators. 
The undeformed osp(2,2) algebra is recovered in the limit
$p=s=r=1$. 
\par Decoupling the fermionic generators from  eq.(\ref{bpm}) leads, 
together with (\ref{b01}), (\ref{b0pm}) and (\ref{b1pm}),  to the
two parameters deformation of gl(2) obtained in ref. \cite{faza}
(a four parameters deformation of gl(2) was further 
obtained in ref.\cite{fanu}). The 
parameter $r$ of the deformation affects only the relations involving 
the fermionic generators.
\section{Representations}
Restricting the three parameters $p,r,s$ in the algebra
above to the case 
\be
    p  = q^{{1\over 2}} \ \ \ , \ \ \ 
    s  = q^{-{1\over 2}} \ \ \ ,  \ \ \ r=1
\ee
leads to a one parameter deformation  of osp(2,2) parametrized by $q$
which we will refer to as osp(2,2)$_q$.
It is equivalent, up to a rescaling of the generators,
to the deformation  discussed in \cite{bgk}.
Like for the Witten type deformation of sl(2),
it is possible to construct some representations 
of osp(2,2)$_q$ in terms
of a finite difference operator $D_q$~:
\be
       D_q f(x) =  { f(x) - f(qx) \over (1-q)x} \ \ \ , \ \ \ 
       D_q x^n = [n]_q x^{n-1} \ \ \ , \ \ \   
       [n]_q \equiv {1 - q^n \over 1-q}  
\label{jack}
\ee 
Some of these representations act on the vector space 
$P(n-1) \oplus P(n)$ (denoting the set of
couples of polynomials whose first (resp. second)
component is of degree at most $n-1$ (resp. $n$) in the variable $x$).
The fermionic generators are represented by
\bea 
     &V_1 = \sigma_-  ,  &V_2 = x \sigma_- \nonumber \\
     &\ov V_1 = q^{-n} (xD_q - [n]_q) \sigma_+  , 
     &\ov V_2 = q^{-1} D_q \sigma_+ 
\label{rep2}
\eea
and the operators representing the  $E_{ij}$'s are easily
constructed through eq.(\ref{bff}). 
The enveloping algebra constructed with the four operators 
(\ref{rep2}) generates all
the   2$\times$2 matrix,  finite difference
operators preserving $P(n-1) \oplus P(n)$. 
Accordingly, these operators are quasi exactly solvable.
Note that more general representations of osp(2,2)$_q$ 
can be constructed \cite{bgk}; they preserve the vector
space  $P(n) \oplus P(n+1) \oplus P(n-1) \oplus P(n)$
(using an obvious notation).
\section{Deformations of osp(1,2)}
Recently a Witten type deformation of the super Lie algebra osp(1,2)
was constructed \cite{chung}.
We showed that this deformation is the only one to fulfil
the three requirements of sect. 2 and we
constructed the representations of 
it this algebra in terms of the operator (\ref{jack}).
Again the space of the representation is 
$P(n-1) \oplus P(n)$ and the two fermionic generators of osp(1,2)
are represented by  2$\times$2 matrix operators.
Using exactly the same notation as in ref.\cite{chung} we find
\be
V_- =
\left(\begin{array}{cc}
0     & D_{q^2}     \\
1     &0  
\end{array}\right)
\ \ \ , \ \ \ \ 
V_+ =
\left(\begin{array}{cc}
0                 &q^{-2n}(x D_{q^2} - [n]_{q^2})     \\
x     &0  
\end{array}\right)
\ee
The three bosonic operators, $H,J_-,J_+$  can then be computed
through the structure of osp(1,2), namely
\be 
    H   =  \{V_-,V_+ \}_q \ \ , \ \ 
    J_- =  \{V_-,V_- \}_q \ \ , \ \
    J_+ =  \{V_+,V_+ \}_q 
\ee
Corresponding to this $2n+1$ dimensional representation, the Casimir operator
(eq.(13) in ref. \cite{chung})
has a value $ C = -{1\over 2} [-n - {1\over 2}]_{q^2}$
\section{Concluding remarks}
The most interesting examples of QES systems are related to
the algebra osp(2,2), e.g. the relativistic Coulomb problem
and the stability of the sphaleron in the abelian Higgs model
\cite{bk}. Moreover, the representations of osp(2,2) formulated in terms
of differential operators provide the building blocks for
the construction of the QES operators preserving the vector space 
$P(m)\oplus  P(n)$ for any  positive integers $m,n$
\cite{bk}. This result and
the recent applications of discrete equations 
in quantum mechanics \cite{wiza}
have encouraged us to study the deformations of osp(2,2)
formulated in terms of quommutators. Our result can be generalized
in several directions; for instance the operators (\ref{rep2}) 
can be used to construct the QES discrete operators preserving
the vector space $P(m)\oplus  P(n)$. The algebraic structure
underlying these operators is probably 
determined by a series of deformed non linear superalgebras
indexed by  $\vert m-n \vert$. Another possibility is to
study the deformations of the Lie superalgebra 
spl(V+1,1) for V$>$1 (remembering the equivalence of osp(2,2) with spl(2,1)).  
The result of ref. \cite{bn} suggests that some
representations will be formulated in terms of  
finite difference operators depending of $V$ variables.
Finally it would  be instructive to relate the two
parameter deformation of osp(2,2) obtained in ref.\cite{para}
to the deformation constructed here.
The occurence of some mapping between the two
structures would allow to transport the Hopf structure 
elaborated in \cite{para} into our deformation.
\vskip 1.5 true cm
\noindent {\bf Acknowledgements.}\\
I gratefully acknowledge  Boucif Abdessalam and Amitabha Chakrabarti 
for discussions and for their hospitality in Ecole Polytechnique
in Paris.

\end{document}